\documentclass{article}
\pdfoutput=1
\usepackage{color, xcolor, colortbl}
\usepackage{graphicx,epstopdf}
\usepackage{geometry}
\usepackage{enumerate}
\usepackage{amsmath,amssymb,amsthm}
\usepackage{algorithm}
\usepackage{algorithmic}
\usepackage{caption}

\usepackage{subcaption}
\usepackage{bm}
\usepackage{appendix}
\usepackage{multirow}
\usepackage{braket}
\usepackage[english]{babel}
\usepackage{hyperref}
\usepackage[capitalize]{cleveref}
\usepackage[sort&compress,square,numbers]{natbib}
\usepackage{adjustbox}
\usepackage{xspace}
\usepackage[roman]{complexity}
\usepackage{soul}
\usepackage{tikz}
\usetikzlibrary{quantikz}
\usepackage{rotating}
\usepackage{setspace}
\usepackage{placeins}
\usepackage{fancyhdr}
\DeclareMathOperator{\Comm}{Comm}

\usepackage{authblk}

\newcommand{\norm}[1]{\left\lVert#1\right\rVert}

\newcommand{\op}{{\rm op}}

\newcommand{\Or}{\mathcal{O}}
\newcommand{\NN}{\mathbb{N}}
\newcommand{\RR}{\mathbb{R}}

\newtheorem{definition}{Definition}

\newtheorem{thm}{\protect\theoremname}
\theoremstyle{plain}
\newtheorem{lemma}[thm]{\protect\lemmaname}
\theoremstyle{plain}
\newtheorem{rem}[thm]{\protect\remarkname}
\theoremstyle{plain}
\theoremstyle{plain}

\theoremstyle{plain}

\providecommand{\definitionname}{Definition}
\providecommand{\assumptionname}{Assumption}
\providecommand{\corollaryname}{Corollary}
\providecommand{\lemmaname}{Lemma}
\providecommand{\propositionname}{Proposition}
\providecommand{\remarkname}{Remark}
\providecommand{\theoremname}{Theorem}

\title{Superconvergence of High-order Magnus Quantum Algorithms}

\author[1,2]{Di Fang}
\author[1,2]{Jiaqi Zhang}

\affil[1]{Department of Mathematics, Duke University}
\affil[2]{Duke Quantum Center, Duke University}

\date{} 

\begin{document}
\maketitle

\begin{abstract}
The Magnus expansion has long been a celebrated subject in numerical analysis, leading to the development of many useful classical integrators. More recently, it has been discovered to be a powerful tool for designing quantum algorithms for Hamiltonian simulation in quantum computing. In particular, surprising superconvergence behavior has been observed for quantum Magnus algorithms applied to the simulation of the Schr\"odinger equation, with the first- and second-order methods exhibiting doubled convergence order. In this work, we provide a rigorous proof that such superconvergence extends to general high-order quantum Magnus algorithms. Specifically, we show that a quantum Magnus algorithm of order $p$ achieves the superconvergence of order $2p$ in time when applying to the Schr\"odinger equation simulation in the interaction picture. Our analysis combines techniques from semiclassical analysis and Weyl calculus, offering a new perspective on the mathematical foundations of quantum algorithms for time-dependent Hamiltonian simulation.
\end{abstract}

\section{Introduction}
Hamiltonian simulation aims to approximate the unitary evolution generated by a Hamiltonian and is central to quantum algorithms for physics, chemistry, and beyond. Of particular interest in both applications and complexity theory is the Schr\"odinger operator
\begin{equation} \label{eq:Sch-eq}
H = -\Delta + V(x),
\end{equation}
which is unbounded and therefore difficult to treat under the standard operator-norm error metric. Indeed, while wave functions lie in $L^2$, the negative Laplacian operator $-\Delta$ maps $H^2$ to $L^2$, but not $L^2$ to $L^2$ and hence is unbounded on $L^2$. This issue of unboundedness cannot be resolved merely by introducing spatial discretization. For illustration, consider the one-dimensional case. After spatial discretization on $N$ grid points, the Laplacian contributes a matrix $A$ whose norm scales like $\|A\|=\mathcal{O}(N^2)$ (e.g., for finite differences), whereas the potential becomes a bounded diagonal matrix $B$. Consequently, if one applies any existing time-independent Hamiltonian simulation algorithm directly to $H=A+B$, the cost inherits a polynomial dependence on $N$ through either $\|H\|$ itself or through commutators such as $\|[A,B]\|=\mathcal{O}(N)$ and $\|[A,[A,B]]\|=\mathcal{O}(N^2)$, even though the per-query circuit depth is only polynomial in the qubit number $n=\log N$. All current time-independent methods (Trotter-type~\cite{ChildsSuTranEtAl2020,ChildsSu2019,WatkinsWiebeRoggeroLee2022,ZhukRobertsonBravyi2023,AftabAnTrivisa2024,Watson2025,WatsonWatkins2024,Somma2015,SahinogluSomma2020,AnFangLin2021,SuHuangCampbell2021,ZhaoZhouShawEtAk2021,ChildsLengEtAl2022,FangTres2023,BornsWeilFang2022,HuangTongFangSu2023,ZengSunJiangZhao2022,GongZhouLi2023,LowSuTongTran2023,ZhaoZhouChilds2024,YuXuZhao2024,ChenXuZhaoYuan2024} and post-Trotter-type~\cite{BerryChildsCleveEtAl2015,KieferovaSchererBerry2019,LowWiebe2019,BerryChildsSuEtAl2020,LowChuang2017,LowChuang2019,GilyenSuLowEtAl2019,ZlokapaSomma2024}) therefore reintroduce a $\mathrm{poly}(N)$ dependence through norm or commutator bounds and fail to achieve the desired $\mathrm{polylog}(N)$ complexity in operator norm error. 
In particular, in our setting all existing time-independent Hamiltonian simulation algorithms necessarily incur a $\mathrm{poly}(N)$ dependence.
Moreover, simulating the Schr\"odinger operator is known to be BQP hard~\cite{ZhengLengLiuWu2024}, which further motivates studying it as a valuable benchmarking problem for assessing the performance and complexity of quantum algorithms.
For our purposes, even if an efficient quantum circuit using $\mathrm{polylog}(N)$ gates per time step is available, we must still prove that the number of time steps does not reintroduce a $\mathrm{poly}(N)$ cost.

To bypass this obstruction, we pass to the interaction picture~\cite{LowWiebe2019}. Conjugating by the kinetic part $A=-\tfrac12\Delta$ replaces the time-independent Hamiltonian by a time-dependent one,
\begin{equation}
H_I(t)=e^{iAt} B e^{-iAt},
\end{equation}
whose operator norm is governed by $\|B\|$ rather than by $\|A\|$ and hence is independent of $N$, if $B$ is bounded. At the same time, the time variation becomes nontrivial: under standard discretizations, $\|\partial_t H_I(t)\|=\|[A,B]\|=\mathcal{O}(N)$. Thus, a fully efficient simulation algorithm requires an efficient time-dependent Hamiltonian simulation solver whose complexity depends only logarithmically on the strength of time variation. For example, this could mean a dependence on $\log(\|\partial_t H(t)\|)$ (or an equivalent measure of time-variation), rather than polynomial dependence. The interaction picture addresses the norm blow-up but leaves large derivatives, so the solver’s derivative dependence must be mild to retain overall $\mathrm{polylog}(N)$ cost.

Recent progress shows that quantum algorithms based on the Magnus expansion meet exactly this need~\cite{AnFangLin2022,FangLiuSarkar2025,FangLiuZhu2025,Borns-WeilFangZhang2025}. These algorithms provide circuits whose gate counts scale $\mathrm{polylog}(N)$ per time step while depending only logarithmically on suitable norms of time derivatives. Beyond this structural advantage, a surprising phenomenon, superconvergence, has been observed for the Schr\"odinger equation in the interaction picture: the first-order quantum Magnus method empirically exhibits second-order accuracy, and the second-order method exhibits fourth-order accuracy~\cite{AnFangLin2022,FangLiuSarkar2025,Borns-WeilFangZhang2025}, with error constants independent of $N$ for smooth bounded potentials. This observation leads to the central question of this paper:

\medskip
\emph{Does superconvergence persist for high-order quantum Magnus algorithms? In particular, is the order gain a uniform doubling from $p$ to $2p$, or does it instead follow a pattern where even orders gain $2$ while odd orders gain $1$?}
\medskip

Establishing a high-order theory faces two distinct challenges. First, high-order Magnus truncations exhibit factorial growth in the number of terms indexed by permutations: the $n$-th term $\Omega_n$ involves sums over $\pi\in S_n$ with coefficients $C_{\pi,n}$ (in the notation used later for our algorithmic description), so the number of distinct terms in the nested-commutator structures scales like $n!$. On the algorithmic side, this factorial explosion has already been eliminated by careful circuit constructions (block-encoding and quantum lookup table-based implementations of truncated Magnus generators), so the circuit side is under control as demonstrated in~\cite{FangLiuZhu2025}. Second, the error analysis in the operator norm must gain the superconvergence cancellation for an arbitrary order of $p$, while the preconstant does not depend sensitively on $\|A\|$ (equivalently, it remains uniform in $N$). This is subtle even at low order and was previously unclear for high order. In particular, one must show that once an efficient quantum circuit using $\mathrm{polylog}(N)$ gates per time step is available, the number of time steps needed to reach a target accuracy does not reintroduce a $\mathrm{poly}(N)$ cost; otherwise the overall algorithm would cease to be $\mathrm{polylog}(N)$.

In this paper, we resolve the above open question positively and completely for the Schr\"odinger equation in the interaction picture. We prove that a quantum Magnus algorithm of order $p$ achieves \emph{superconvergence of order $2p$ in time}. Precisely, letting
\begin{equation} \label{eq:U}
U(T,0) \;=\; \mathcal{T}\exp\!\left(-i \int_{0}^{T} H(s)\,ds\right)
\end{equation}
be the exact evolution operator, and
\begin{equation}\label{eq:U_p}
U_p(T,0) \;=\; \prod_{k=0}^{L-1}\exp\!\left(\sum_{n=1}^{p}\Omega_n(t_{k+1},t_k)\right)
\end{equation}
be the global $p$-th order Magnus approximation built from step size $h=T/L$ (see \cref{sec:algo_overview}), we show that for $A=-\tfrac12\Delta$ and $B=V(x)$ with $V(x)\in S(1)$ (bounded together with all derivatives) the error bound
\begin{equation}
\|U(T,0)-U_p(T,0)\|\;\le\; C_{V,p}\,T\,h^{2p}
\end{equation}
holds, where the constant $C_{V,p}$ depends only on $V$ and finitely many of its derivatives, but is independent of norm of $A$. This implies with the algorithm discretization, the error is bounded independent of the spatial discretization $N$.

Our analysis develops a uniform commutator calculus for the interaction-picture Hamiltonian. 
The starting point is a pull-out identity that conjugates all time labels to a common
reference time. Concretely, any nested commutator built from
$H_I(t)=e^{iAt}Be^{-iAt}$ at times $t_1,\ldots,t_q$ can be rewritten, up to an overall
conjugation by $e^{iAs_0}$, as a nested commutator among the shifted operators
$B_{t_j-s_0}=e^{iA(t_j-s_0)}Be^{-iA(t_j-s_0)}$ with innermost entry $B_{s_0}$.
In particular, taking $s_0=0$ shows that such commutators reduce to combinations of
$B_0=B$ together with $B_{t_j}$.
We then work semiclassically with Weyl quantization $\op_h$ on $T^*\RR^d$ and exploit two basic facts recalled in \cref{sect:prelim}: the exact Egorov theorem for the quadratic kinetic symbol $a(p)=p^2/2$, which gives $B_u=\op_h\big(v(x-u p)\big)$ exactly, and the semiclassical commutator rule $[\op_h(a),\op_h(b)]=\tfrac{h}{i}\op_h(\{a,b\})+h^2\op_h(S(1))$. Iterating these tools shows that any $q$-layer commutator among the $B_u$ has the form $h^{q}\op_h(S(1))$ with uniform symbol bounds. An application of Calder\'on-Vaillancourt Theorem then yields operator-norm control that is uniform in $h$. At the level of Magnus error representations this uniform $h$-gain doubles the formal order from $p$ to $2p$ and, crucially, the constants depend only on $V$ and its derivatives, not on $A$ or on $N$. The result gives a rigorous explanation of the superconvergence phenomenon observed at low orders and extends it to all orders.

To connect this accuracy theory to overall complexity, note that the interaction picture yields $\|H_I(t)\|=\|B\|$, while the large contributions come from time variation. Modern time-dependent solvers based on Magnus expansions can be arranged so that their dependence on such time-variation measures enters only logarithmically. When this is combined with our $2p$-order error bound, the step count $L=T/h$ required to meet a target tolerance depends only on $T$ and $V$ (through $C_{V,p}$), and not on $N$. This completes the complexity picture: having $\mathrm{polylog}(N)$ per-step cost and an $N$-independent step count ensures that the overall gate complexity is $\mathrm{polylog}(N)$.

Finally, we remark that the Magnus expansion is not only a basis for algorithms but also an analytical and theoretical tool in classical numerical analysis and other areas~\cite{HochbruckLubich2003,Thalhammer2006,BlanesCasasThalhammer2017,AnFangLin2022,FangLiuSarkar2025,BlanesCasasOteoRos2009,SharmaTran2024,BosseChildsEtAl2024,CasaresZiniArrazola2024,ApelCubittOnorati2025}. Our estimates are formulated directly at the level of symbols and commutators, rather than through Taylor series-based argument. 
This perspective highlights structural cancellations, and the resulting techniques may be of independent interest for classical Magnus integrators as well as for other applications where commutator effects and time-dependent structure are essential.

The remainder of the paper is organized as follows. \cref{sect:prelim} collects the semiclassical notation and tools we use, including symbol classes, Weyl quantization, the exact Egorov theorem, and the Calder\'on-Vaillancourt Theorem. We then review the $p$-th order quantum Magnus algorithm and its circuit implementation, focusing on the mathematical form rather than circuit engineering details. Next, we formulate the Schr\"odinger equation in the interaction picture and state the commutator estimates that underlie our analysis. The main superconvergence theorem and its proof appear in the section that follows, where we derive $2p$-th order bounds with constants depending only on $V$ and its derivatives. 
We conclude with numerical illustrations and remarks about extensions, including the discrete setting and further applications where the Magnus series can be used as an analytical tool of independent theoretical interest.

\section{Preliminary} \label{sect:prelim}
In this section, we introduce several notational conventions and fundamental properties that will be used throughout the paper. These preliminaries provide a common framework for our arguments and serve as a reference point for subsequent proofs. By laying out these tools in advance, we ensure that the subsequent analysis can proceed in a clear and consistent manner.

Before introducing the symbol class $S(m)$, we revisit a few notational conventions that we will use throughout. We write $z=(x,p)\in\RR^{2d}$ and use the Japanese brackets
\begin{equation}
\langle u\rangle := (1+|u|^2)^{1/2}\quad (u\in\RR^k).
\end{equation}
Now, let $m$ be a order parameter in the usual definition of symbol class \cite{zworski2022semiclassical}, which means that there exist constants $C>0$ and $N\ge 0$ such that
\begin{equation}\label{eq:order}
m(w)\le C\,\langle z-w\rangle^{N}\,m(z)\qquad\text{for all }w,z\in\RR^{2d}.
\end{equation}
A particularly useful example is \cite{BornsWeilFang2022}
\begin{equation}
m(x,p)=\langle x\rangle^{a}\,\langle p\rangle^{b}
\quad\text{for any }a,b\in\RR,
\end{equation}
where the Japanese brackets are defined as above.

\begin{definition}[Symbol class $S(m)$]\label{def:S1}
Phase space is $T^*\mathbb{R}^d$ with coordinates $(x,p)$ (position $x\in\mathbb{R}^d$, momentum $p\in\mathbb{R}^d$; thus $T^*\RR^d\cong\RR^d_x\times\RR^d_p$).
We set
\begin{equation}
S(m):=\Big\{a\in C^\infty(T^*\RR^d): \text{ for all multi-indices } \alpha,\beta, \exists\  C_{\alpha\beta} \ \text{s.t.} \ \big|\partial_x^\alpha\partial_p^\beta a(x,p)\big|\le C_{\alpha\beta} m(x,p) \Big\}.
\end{equation}
\end{definition}
In particular, $S(1)$ means that all mixed derivatives of $a$ are bounded on $T^*\RR^d$. When we choose $\langle x\rangle^a\langle p\rangle^b$, the symbol class can include polynomials in $x$ and $p$.

\begin{definition}[Weyl quantization]\label{def:Weyl}
For $a\in S(m)$ and $h\in(0,1]$, its Weyl quantization is
\begin{equation}
\op_h(a)u(x)
:= \frac{1}{(2\pi h)^d}\int_{\RR^d}\!\!\int_{\RR^d}
e^{\frac{i}{h}(x-y)\cdot p}\,
a\!\left(\frac{x+y}{2},p\right) u(y)\,dp\,dy,\qquad u\in \mathcal{S}(\RR^d).
\end{equation}
Here, $\mathcal{S}(\RR^d)$ denotes the Schwartz space of smooth functions that decay rapidly at infinity together with all their derivatives.  
\end{definition}
More generally, the definition can extend to $a\in\mathcal{S}'(\RR^{d})$, the space of tempered distributions, by duality, so that $\op_h(a)$ defines a continuous map on $\mathcal{S}(\RR^d)$. We review the following helpful properties. For more details of symbol class and Weyl quantization, see e.g.,~\cite{zworski2022semiclassical,Martinez2002}.

\begin{lemma}[Semiclassical commutator rule]\label{lem:comm_rule}
For the commutator $[\op_h(a),\op_h(b)]$ with $a,b\in S(1)$, there exists a symbol $c\in S(1)$ such that
\begin{equation}\label{eq:semi-comm}
\big[\op_h(a),\op_h(b)\big] \;=\; \frac{h}{i}\,\op_h\!\big(\{a,b\}\big)\;+\;h^2\,\op_h(c).
\end{equation}
Here $\{a,b\}$ is the Poisson bracket,
\begin{equation}
\{a,b\}\;:=\;\sum_{j=1}^d\!\left(\frac{\partial a}{\partial p_j}\frac{\partial b}{\partial x_j}
-\frac{\partial a}{\partial x_j}\frac{\partial b}{\partial p_j}\right),
\end{equation}
which also lies in $S(1)$ whenever $a,b\in S(1)$.
\end{lemma}

\begin{lemma}[Exact Egorov Theorem]\label{lem:egorov}
For any $v\in S(1)$ depending only on $x$ and any $s\in\RR$,
\begin{equation}
e^{ish\Delta}\op_h(v(x))e^{-ish\Delta} =\op_h\big(v(x-sp)\big),
\end{equation}
where $v(x-sp)\in S(1)$ and $\Delta$ is the Laplacian operator.
\end{lemma}
In particular, setting $V(x):=\op_h(v(x))$, we see that the time-evolved observable in our case is
\begin{equation}
V_A(t):=e^{iAt}Ve^{-iAt}.
\end{equation}
By the above lemma, for $t=hs$ this becomes
\begin{equation}
V_A(hs)=\op_h\big(v(x-sp)\big),
\end{equation}
where the symbol $v(x-sp)\in S(1)$ uniformly. As we will see in later discussion, the exact Egorov Theorem thus describes the behavior of $V$ under the time evolution, and this is directly connected to the interaction picture, which we will examine in detail in a later section.

\begin{thm}[Calder\'on-Vaillancourt Theorem]\label{thm:CV}
Let $n\ge 1$ and let $a\in S(1)$ be a symbol on $T^*\RR^n$.
Then the Weyl operator $\op_h(a):L^2(\RR^n)\to L^2(\RR^n)$ is a bounded (continuous linear) map. 
More precisely, there exist constants $M_n$ and $C_n$, depending only on $n$, such that
\begin{equation}
\big\|\op_h(a)\big\|_{L^2\to L^2}
\;\le\; C_n \Bigg(\sum_{|\alpha|\le M_n}\ 
\sup_{(x,p)\in T^*\RR^n} \big|\partial_{x,p}^{\alpha} a(x,p)\big|\Bigg).
\end{equation}
Here $\partial_{x,p}^{\alpha}$ denotes any mixed derivative in the variables $(x,p)$.
The constants are independent of $h$.
\end{thm}

\section{Algorithmic Overview}  \label{sec:algo_overview}

In this section, we first review the basics of Magnus expansion, and then revisit the $p$-th order quantum Magnus algorithm introduced in~\cite{FangLiuZhu2025}. Rather than reproducing the circuit constructions, we focus on the mathematical formulation of the algorithm, which serves as the foundation for our subsequent analysis of superconvergence.

Let $U(t)$ solve $\dot U(t)=A(t)U(t)$ with $U(0)=I$, and set $A(t):=-iH(t)$.
The time-ordered exponential admits the Magnus representation
\begin{equation}
\mathcal{T}e^{\int_{0}^{t}A(s)\,ds}=e^{\Omega(t)},\qquad
\Omega(t)=\sum_{n=1}^{\infty}\Omega_n(t)=:\Omega_{(\infty)}(t),
\end{equation}
where $\mathcal{T}$ denotes the time-ordering operator, and the first three terms are
\begin{align*}
\Omega_1(t) &= \int_{0}^{t} A(t_1)\,dt_1,\\
\Omega_2(t) &= \frac{1}{2}\!\int_{0}^{t}\!dt_1\!\int_{0}^{t_1}\!dt_2\,[A(t_1),A(t_2)],\\
\Omega_3(t) &= \frac{1}{6}\!\int_{0}^{t}\!dt_1\!\int_{0}^{t_1}\!dt_2\!\int_{0}^{t_2}\!dt_3\,
\big([A(t_1),[A(t_2),A(t_3)]]+[A(t_3),[A(t_2),A(t_1)]]\big),
\end{align*}
with $[X,Y]=XY-YX$.
We truncate after $p$ terms and write
\begin{equation}
\Omega_{(p)}(t):=\sum_{n=1}^{p}\Omega_n(t),\qquad
U_p(t):=\exp\big(\Omega_{(p)}(t)\big).
\end{equation}
This is the $p$-th order Magnus expansion, which has wide applications in classical algorithm design and as an analytical tool for studying physical and chemical systems (see, e.g., the book~\cite{BlanesCasas2017book} and the reviews~\cite{BlanesCasasOteoRos2009,HairerHochbruckIserlesLubich2006}).

We now review the high-order quantum Magnus algorithm as proposed in ~\cite{FangLiuZhu2025}. On a subinterval $[t_k,t_{k+1}]$ of length $h$, we use the same formulas with the
integration limits replaced by $t_k$ and $t_{k+1}$, and denote
\begin{equation}
\Omega_{(p)}(t_{k+1},t_k):=\sum_{n=1}^{p}\Omega_n(t_{k+1},t_k),\qquad
U_p(t_{k+1},t_k):=\exp\!\big(\Omega_{(p)}(t_{k+1},t_k)\big).
\end{equation}
In each timestep we approximate the exact local dynamics by the $p$th-order Magnus unitary $U_p(t_{k+1},t_k)$; the full evolution is obtained by composing these per-step unitaries. 
We design a quantum algorithm using the high-order Magnus expansion at each of the short time interval $[t_j, t_{j}+h]$, namely,
\begin{equation}\label{eq:U_exact_apprx_U_p}
        U(t_j+h,t_j) 
    = \mathcal{T}\exp\!\left(-i \int_{t_j}^{t_j+h} H(s)\, ds\right) 
    \;\;\approx\;\; 
    \underbrace{\exp\!\left(\Omega_{(p)}(t_j+h, t_j)\right)}_\text{Magnus expansion}
    \approx \underbrace{\exp\!\left(\tilde{\Omega}_{(p)}(t_j+h, t_j)\right)}_\text{using numerical quadrature},
\end{equation}
where the $p$-th Magnus expansion is given by
$
  \Omega_{(p)}(t_j+h, t_j) =  \sum_{n=1}^p \Omega_n(t_j+h, t_j) 
$
and
\begin{equation*}
    \Omega_n(t_j+h,t_j) = \sum_{\pi\in S_n} C_{\pi,n} \int_{t_j}^{t_j+h} ds_1 \int_{t_j}^{t_j+s_1} ds_2 \cdots \int_{t_j}^{t_j+s_{n-1}} ds_{n} A(s_{\pi(1)})A(s_{\pi(2)})\cdots A(s_{\pi(n)}),
\end{equation*}
where $$ C _{\pi,n}= \frac{(-1)^{d_a(\pi)}}{n} \times \frac{1}{ {n-1 \choose d_a(\pi)}}, \quad d_a(\pi)=\left|\{\,i\in\{1,\dots,n-1\}|\pi(i)>\pi(i+1)\}\right|.$$ Here $d_a(\pi)$ denotes the number of descents of permutation $\pi$ and $n$ denotes the length of permutation.

We choose a sufficiently large number of numerical quadrature points (denoted by $M$) so that the quadrature error is comparable to the local truncation error from truncating the Magnus series. The quantum circuit described in~\cite{FangLiuZhu2025} incurs only logarithmic cost in $M$. In other words, we can safely increase the number of quadrature points without introducing significant quantum-cost overhead. Therefore, in what follows, we focus on estimating the error arising from the Magnus expansion itself.
The algorithm is described below.

\begin{algorithm}[H]
\caption{Revisit of the $p$-th order quantum Magnus algorithm}
\begin{itemize}
\item \textbf{Inputs}: Final simulation time $T>0$; number of steps $L\in\mathbb{N}$ (with step size $h = T/L$); Magnus order $p$; time-dependent Hamiltonian $H(t)$.
\item \textbf{Initialization}: Partition the interval $[0,T]$ into $L$ subintervals with time grid points $t_k = k h$ for $0 \leq k \leq L$. The $k$-th time step corresponds to the subinterval $[t_k, t_{k+1}]$. 
\item For $k=0,1,\dots,L-1$ (time step $[t_k,t_{k+1}]$):

Construct the p-th order Magnus expansion $\Omega_{(p)}(t_{k+1}, t_k)$ using a large number of quadrature points using a linear combination of unitary (LCU)-type of quantum circuit (the quantum cost is only $\log(M) = \log(\norm{H'(t)})$), and then implement the approximate Magnus unitary \begin{equation}\label{eq:def_Up_tk_tk+1}
    U_p(t_{k+1}, t_k)) = \exp(\Omega_{(p)}(t_{k+1}, t_k))
\end{equation} using a Quantum Singular Value Transformation (QSVT) circuit. All circuits are explicitly constructed and presented in~\cite{FangLiuZhu2025}. 

\item \textbf{Output}: The final approximate evolution operator is the product of the per-step Magnus unitaries:
$$
U_p(T,0)=\prod_{k=0}^{L-1}U_p(t_{k+1}, t_k))   \;=\; \prod_{k=0}^{L-1} \exp\!\big(\,\Omega_{(p)}(t_{k+1}, t_k)\big).
$$
\end{itemize}
\end{algorithm}

While a general complexity analysis for arbitrary time-dependent Hamiltonians $H(t)$ was given in~\cite{FangLiuZhu2025} where we review \cref{thm:magnus_trunction_finite_comm}. The algorithm achieves pth-order accuracy, this work focuses on the case where the underlying problem is the Schr\"odinger PDE in \cref{eq:Sch-eq} %
in the interaction picture. In this setting, we establish rigorous estimates for the required step size (equivalently, the number of steps $L$), and prove that the $p$-th order quantum Magnus algorithm attains $2p$-th order of convergence. Moreover, we show that the error prefactor depends only on the potential $V(x)$ and its derivatives.

\begin{lemma}[Local Truncation Error of $p$-th Order Magnus Expansion \cite{FangLiuZhu2025}] \label{thm:magnus_trunction_finite_comm}
 Let $U(t_{j+1},t_j)$
be the exact short-time propagator over the time interval $[t_j, t_{j+1}]$ given by \cref{eq:U_exact_apprx_U_p}, and let $U_p(t_{j+1},t_j)$ denote the approximation obtained by applying the $p$-th order Magnus expansion given by \cref{eq:def_Up_tk_tk+1} with time step size $h = T/L$
and $t_j = j h$.
Let $\alpha_{\mathrm{comm},q}$ denote the maximum norm (taken over time $[t_j, t_{j+1}]$ and over all nested commutators of grade $q$) of the operator $H(t)$. That is,
 \begin{equation}\label{eq:def_alpha_comm_q}
\alpha_{\mathrm{comm},q} := \sup_{\tau_1, \ldots, \tau_q \in [t_j,t_{j+1}]} \max_{C \in \mathcal{C}_q} \big\| C\big(H(\tau_1), \ldots, H(\tau_q)\big) \big\|,
 \end{equation}
where $\mathcal{C}_q$ denotes the set of all nested commutators of grade $q$. 
Then the local truncation error is bounded by
\begin{equation} \label{eq: LTE bound}
  \norm{U(t_{j+1},t_j) - U_p(t_{j+1},t_j)} \leq  C  \sum_{q = p+1}^{p^2+2p} \alpha_{\mathrm{comm},q}  h^{q} ,
\end{equation}
where $C$ is some constant depending only on $p$.
\end{lemma}

\section{Schr\"odinger Equation in the Interaction Picture}\label{subsec:IP}
In this section, we set up the interaction picture framework for the Schr\"odinger equation.
Our goal is to decompose the Hamiltonian into kinetic and potential parts, introduce the
interaction-picture Hamiltonian, and establish the equivalence between the Schr\"odinger
evolution and its interaction-picture formulation. These constructions will provide the
foundation for the Magnus expansion analysis developed in the subsequent sections.
We decompose the Hamiltonian into kinetic and potential parts
\begin{equation}
  A = -\tfrac{1}{2} \Delta, \qquad B = V(x), \qquad H = A + B .
\end{equation}
The Schr\"odinger evolution operator $U(t)$ satisfies
\begin{equation}\label{eq:sch-AB}
  i \partial_t U(t) = (A+B)\, U(t), \qquad U(0) = I .
\end{equation}
To define the interaction picture, we introduce
\begin{equation}\label{eq:HI-def}
  H(t):= H_I(t) = e^{iAt}\, B \, e^{-iAt},
\end{equation}
and let $U_I(t)$ be the solution of
\begin{equation}\label{eq:UI-def}
  i \partial_t U_I(t) = H_I(t)\, U_I(t), \qquad U_I(0) = I .
\end{equation}
The two pictures are equivalent: setting
\begin{equation}\label{eq:U-UI}
  U(t) = e^{-iAt}\, U_I(t)
\end{equation}
yields \cref{eq:sch-AB} from \eqref{eq:UI-def}, and conversely 
\begin{equation} \label{eq:U_IP}
U_I(t)=e^{iAt} U(t)
\end{equation}
satisfies
\cref{eq:UI-def} whenever $U(t)$ satisfies \cref{eq:sch-AB}. Thus the interaction picture is obtained
by conjugating the potential part $B$ under the unitary evolution generated by $A$. For our case, the interaction picture is well-defined as proved in~\cite[Section 4.1]{FangLiuSarkar2025}.

For later reference, we also record the common notation
\begin{equation}\label{eq:Bu}
  B_u^A = e^{iAu} B e^{-iAu}, \qquad u \in \mathbb R ,
\end{equation}
and in particular
\begin{equation}\label{eq:BA}
  B_{t}^A = e^{iAt} B e^{-iAt} .
\end{equation}
With this notation, $H_I(t) = B_{t}^A$, so the interaction-picture Hamiltonian is the $A$-evolved
potential.

\section{Proof of Superconvergence}

Let $\alpha_{\mathrm{comm},q}$ denote the maximum norm (taken over time and over all nested commutators of grade $q$) of the operator $A(t)$. We want to show:
 \begin{equation}
 \sup_{\tau_1, \ldots, \tau_q \in [t_j,t_{j+1}]} \max_{C \in \mathcal{C}_q} \big\| C\big(H(\tau_1), \ldots, H(\tau_q)\big) \big\| \leq C_V h^{q-1},
 \end{equation}
where $\mathcal{C}_q$ denotes the set of all nested commutators of grade $q$. 

\begin{equation}
    H(t) = e^{iAt} B e^{-iAt}, \quad A = -\frac{1}{2}\Delta, \quad B = V(x) \in S(1).
\end{equation}

\begin{definition}[Grade-$n$ commutators w.r.t.\ $B$] We first define grade-$n$ (gd) as the nested commutators that have $n$ occurrences of evolved $B$ terms; equivalently it has $n-1$ commutator layers.
\end{definition}
Let $A$ and $B$ be time-independent operators. For $u\in\mathbb{R}$, we let
\begin{equation}
B_u^{A}=e^{iAu}Be^{-iAu}.
\end{equation}
For notational simplicity, we drop $A: B_u \triangleq B_u^A$,  and only specify this when it is not evolved under $A$.
In particular, $H_I(t)=e^{iAt}Be^{-iAt}=B_t$. Also, as $B=V(x)$ in our case, we sometimes also use $V_t$ for $H_I(t)$.

For the rest of the paper, we adopt the following notation convention:
\begin{definition}
For $n\ge 1$, define the grade-$n$ left-normed commutators (with respect to $B$) recursively by
\begin{equation}
\Comm^{B,\mathrm{gd}}_{1}(s_0):=B_{s_0},\qquad
\Comm^{B,\mathrm{gd}}_{n+1}(s_0;\,s_1,\dots,s_n)
:=\big[\,B_{s_n},\,\Comm^{B,\mathrm{gd}}_{n}(s_0;\,s_1,\dots,s_{n-1})\,\big].
\end{equation}
\end{definition}
Equivalently,
\begin{equation}
\Comm^{B,\mathrm{gd}}_{n}(s_0;\,s_1,\dots,s_{n-1})
= \operatorname{ad}_{B_{s_{n-1}}}\cdots\operatorname{ad}_{B_{s_1}}\!\big(B_{s_0}\big).
\end{equation}

If the innermost term is $B_{s_0}$ with $s_0=0$ (so $B_{s_0}=B$ is time-independent), we drop this argument and, for $q\ge 0$, write
\begin{equation}
\Comm_{q+1}(s_1,\dots,s_q):=\Comm^{B,\mathrm{gd}}_{q+1}(0;\,s_1,\dots,s_q).
\end{equation}
Unless otherwise specified, we abbreviate superscripts and simply write $\Comm_q$.
With this convention a grade-$(q+1)$ object corresponds to a $q$-layer commutator. In particular,
\begin{equation} \label{eq:def_Comm}
\Comm_2(s_1)=[B_{s_1},B],\quad
\Comm_{k+2}(s_1,\dots,s_{k+1})=\big[B_{s_{k+1}},\,\Comm_{k+1}(s_1,\dots,s_k)\big].
\end{equation}

\begin{lemma}[Pull-out rule for $q$-layer commutators]\label{lem:pullout}
Let $A$ and $B$ be general time-independent operators and set
\begin{equation}
    H(t)=e^{iAt}Be^{-iAt}.
\end{equation}
Then for times $t_1,\dots,t_q,t$,
\begin{equation}\label{eq:pullout}
\operatorname{ad}_{H(t_{q})}\cdots\operatorname{ad}_{H(t_1)}\big(H(t)\big)
= e^{\,iAt}\,\Comm_{q+1}(t_1-t,\dots,t_q-t)\,e^{-\,iAt},
\end{equation}
where the left-normed commutators are defined recursively as in \cref{eq:def_Comm}.
Equivalently,
\begin{equation}
\Comm_{q+1}(s_1,\dots,s_q)=\operatorname{ad}_{B_{s_q}}\cdots\operatorname{ad}_{B_{s_1}}(B).
\end{equation}
\end{lemma}

\begin{proof}
We proceed by induction on the \emph{grade} $q\ge 2$ (recall: grade $q$ corresponds to $q-1$ layers).
We use the conjugation identity
\begin{equation}\label{eq:push}
[X,UZU^{-1}] = U[U^{-1} X U,Z]U^{-1},
\end{equation}
valid for any invertible $U$ and bounded operators $X,Z$. This follows by a one-line expansion:
\begin{equation}
[X,UZU^{-1}]
= XUZU^{-1}-UZU^{-1}X
= U\big(U^{-1}XUZ - ZU^{-1}XU\big)U^{-1}.
\end{equation}
\emph{Base case (grade $2$; i.e.\ one layer).}
Using $H(t)=e^{\,iAt}Be^{-\,iAt}$ and \cref{eq:push},
\begin{equation}
\operatorname{ad}_{H(t_1)}\big(H(t)\big)
= [H(t_1),\,e^{\,iAt}Be^{-\,iAt}]
= e^{\,iAt}\big[\,e^{-\,iAt}H(t_1)e^{\,iAt},B\big]e^{-\,iAt}.
\end{equation}
Since $H(t_1)=B_{t_1}$ and $e^{-\,iAt}B_{t_1}e^{\,iAt}=B_{t_1-t}$,
\begin{equation}
\operatorname{ad}_{H(t_1)}\big(H(t)\big)
= e^{\,iAt}[B_{t_1-t},B]e^{-\,iAt}
= e^{\,iAt}\,\Comm_2(t_1-t)\,e^{-\,iAt},
\end{equation}
which is \cref{eq:pullout} with $q=2$. \\
\emph{Inductive step.}
Assume \cref{eq:pullout} holds for some grade $q\ge 2$, i.e.,
\begin{equation}
\operatorname{ad}_{H(t_{q-1})}\cdots\operatorname{ad}_{H(t_1)}\big(H(t)\big)
= e^{\,iAt}\,\Comm_q(t_1-t,\dots,t_{q-1}-t)\,e^{-\,iAt}.
\end{equation}
Then, using \cref{eq:push} and $e^{-\,iAt}H(t_q)e^{\,iAt}=B_{t_q-t}$,
\begin{equation}
\begin{aligned}
\operatorname{ad}_{H(t_q)}\!\Big(
\operatorname{ad}_{H(t_{q-1})}\cdots\operatorname{ad}_{H(t_1)}(H(t))\Big)
&= \big[A(t_q),\,e^{\,iAt}\Comm_q e^{-\,iAt}\big]\\
&= e^{\,iAt}\,\big[\,e^{-\,iAt}H(t_q)e^{\,iAt},\,\Comm_q\,\big]\,e^{-\,iAt}\\
&= e^{\,iAt}\,\big[\,B_{t_q-t},\,\Comm_q\,\big]\,e^{-\,iAt}\\
&= e^{\,iAt}\,\Comm_{q+1}(t_1-t,\dots,t_q-t)\,e^{-\,iAt},
\end{aligned}
\end{equation}
which is \cref{eq:pullout} with the grade advanced from $q$ to $q+1$.
\end{proof}

Starting from this lemma, we work with the specific operators 
\begin{equation}
A = -\frac{1}{2}\Delta, \quad B = V(x).
\end{equation}
Thus, all commutators are taken with respect to $V$, so we follow the same notational simplification and write 
\(\Comm_n^V := \Comm_n\). From \cref{lem:left-normed-cont} and \cref{thm:main-L2} onward, we fix $B=V$ and use the two interchangeably, 
in particular recall $B_u = V_u$ for all $u$.

\begin{lemma}[Left-normed $n$-grade ($n\!-\!1$ layers) commutator]\label{lem:left-normed-cont}
Let $0<h\le 1$, let $A=-\frac{h^2}{2}\Delta$, and let $V(x)\in S(1)$.
For a rescaled time $\tilde{s}_j\in[0,1]$ and $s_j=h\tilde{s}_j$ for all $j\in\NN$.
Then, for every $n\ge1$, there exists $b_{\tilde s_1,\ldots,\tilde s_n}(x,p)\in S(1)$ such that
\begin{equation}\label{eq:left-normed-form}
\Comm_{n+1}(h\tilde s_1,\ldots,h\tilde s_n)
= h^{\,n}\,\op_h\!\big(b_{\tilde s_1,\ldots,\tilde s_n}(x,p)\big).
\end{equation}
where $\Comm_{n+1}(\cdot)$ is defined as in \cref{eq:def_Comm}.

\begin{proof}
First notice that $V(x)=\op_h(V(x))$ for $V(x)\in S(1)$. By the Exact Egorov theorem in \cref{lem:egorov} for the quadratic symbol $a(p)=p^2/2$ and the evolution
notation recalled above, since $\op_h(a(p))=-\frac{h^2}{2}\Delta = h^2 A$ with $a$ quadratic, 
\begin{equation}
V_{h\tilde s}=e^{\,iAh\tilde s}\,V\,e^{-\,iAh\tilde s}
= e^{\frac{i}{h}\tilde s\,\op_h(a(p))}\,\op_h(v(x))\,e^{-\frac{i}{h}\tilde s\,\op_h(a(p))} 
= \op_h\!\big(v(x-\tilde s p)\big)=\op_h(w_{\tilde s}).
\end{equation}
The map $(x,p)\mapsto(x-\tilde s p,p)$ preserves $S(1)$ uniformly for $0\le\tilde s\le1$, hence
$w_{\tilde s}\in S(1)$.

\textbf{Grade 2.}
Using the semiclassical commutator rule as in \cref{lem:comm_rule},
\begin{equation}
[\op_h(f),\op_h(g)]=ih\,\op_h(\{f,g\})+h^2\op_h(r_{f,g}),\qquad r_{f,g}\in S(1),
\end{equation}
Hence, for any $0\le \tilde s_1,\tilde s_2\le1$,
\begin{equation}
\{w_{\tilde s_1},V\}
=\sum_{j=1}^d\!\big(\partial_{p_j}w_{\tilde s_1}\,\partial_{x_j}V
-\partial_{x_j}w_{\tilde s_1}\,\partial_{p_j}V\big)
=-\tilde s_1\,\nabla V(x-p\tilde s_1)\!\cdot\!\nabla V(x)\in S(1),
\end{equation}
Therefore,
\begin{equation}
\Comm_2(h\tilde s_1)
=\big[V_{h\tilde s_1}, V\big]
= h\,\op_h\!\Big(i\,\{w_{\tilde s_1},V\}+h\,r_1\Big),
\end{equation}
which is $h^1$ times Weyl quantization of an $S(1)$ symbol, and thus matches \cref{eq:left-normed-form} for $n=1$.

\textbf{General grade $n\ge 2$.}
Assume \cref{eq:left-normed-form} holds for some $n$ with a symbol 
$b_{\tilde s_1,\ldots,\tilde s_n}(x,p)\in S(1)$.
Then the uniform symbol bounds and the recursive choice may be written as
\begin{equation}\label{eq:uniform-S1-moved}
\sup_{\substack{0\le \tilde s_1,\dots,\tilde s_n\le 1 \\ (x,p)\in \mathbb{R}^{2d}}}
\big|\partial_{x,p}^{\alpha} b_{\tilde s_1,\ldots,\tilde s_n}(x,p)\big|
\le C_{n,\alpha},
\end{equation}
and we can set 
\begin{equation}\label{eq:recursion-bn-moved}
b_{\tilde s_1,\ldots,\tilde s_{n+1}}(x,p)
= i\,\{w_{\tilde s_{n+1}},\,b_{\tilde s_1,\ldots,\tilde s_n}(x,p)\}
\;+\;h\,r_{n+1}(\tilde s_1,\ldots,\tilde s_{n+1}),
\qquad r_{n+1}\in S(1),
\end{equation}
with the same type of bounds as in \cref{eq:uniform-S1-moved}.  
Finally, applying the semiclassical commutator rule in \cref{eq:semi-comm} gives
\begin{equation}
\begin{aligned}
\Comm_{n+2}(h\tilde s_1,\ldots,h\tilde s_{n+1})
&=\big[V_{h\tilde s_{n+1}},\,\Comm_{n+1}(h\tilde s_1,\ldots,h\tilde s_{n})\big]\\
&=\big[\op_h(w_{\tilde s_{n+1}}),\,h^{n}\op_h(b_{\tilde s_1,\ldots,\tilde s_n}(x,p))\big]\\
&=h^{n}\!\left(ih\,\op_h(\{w_{\tilde s_{n+1}},b_{\tilde s_1,\ldots,\tilde s_n}(x,p)\})
+h^2\op_h(r_{n+1})\right) \\
&=h^{n+1}\,\op_h\!\Big(i\{w_{\tilde s_{n+1}},b_{\tilde s_1,\ldots,\tilde s_n}(x,p)\}+h\,r_{n+1}\Big),
\end{aligned}
\end{equation}
which is \cref{eq:left-normed-form} with $n\mapsto n+1$. This completes the proof.
\end{proof}
\end{lemma}

\begin{thm}[Main $L^2\to L^2$ norm estimate for grade $p$ commutator]\label{thm:main-L2}
Let $0<h\le 1$, and let $H(\tau)$ defined as \cref{eq:HI-def}.
Let $\mathcal{C}_q$ denotes the set of all nested commutators of grade $q$. In particular, fix $p\ge1$, times $\tau_1,\ldots,\tau_q\in[0,h]$, and any element in $\mathcal{C}_q$ is a nested commutator
$C\big(H(\tau_1),\ldots,H(\tau_p)\big)$ built from the $p$ time-labelled occurrences
$H(\tau_1),\dots,H(\tau_p)$ (in an arbitrary bracketing). Thus $C$ has exactly $p-1$ layers.
Then there exists a constant $C_{V,p}>0$ such that %
\begin{equation}\label{eq:main-L2}
\sup_{\tau_1,\ldots,\tau_q\in[0,h]}\ \max_{C\in \mathcal{C}_p}\
\big\|\,C\!\big(H(\tau_1),\ldots,H(\tau_p)\big)\,\big\|_{L^2(\mathbb{R}^d) \to L^2(\mathbb{R}^d)} \le\ C_{V,p}\,h^{\,p-1}.
\end{equation}
where $C_{V, p}$ is a constant depending only on $V(x)$ and its derivatives, $p$ and the dimension $d$, and uniformly in $h$.
\end{thm}
\begin{rem}[Equivalence of layers and powers of $h$]
A $q$-layer nested commutator built from the time-labelled operators
$\{H(\tau_j)\}_{j=1}^{q+1}$ means that exactly $q$ commutator brackets are applied among
$q+1$ occurrences of $H(\tau_j)$ (with arbitrary bracketing). Equivalently, such a commutator
is of grade $p=q+1$. We can denote the set of all $q$-layer nested commutators as $\tilde{\mathcal{C}}_q$, which is equivalent to $\mathcal{C}_{q+1}$. Estimate \cref{eq:main-L2} with $p=q+1$ then reads
\begin{equation}
\sup_{\tau_1,\ldots,\tau_q\in[0,h]}
\ \max_{C\in \tilde{\mathcal{C}}_q}\ 
\bigl\|C\!\big(H(\tau_1),\ldots,H(\tau_{q+1})\big)\bigr\|_{L^2\to L^2}
\ \le\ C_{V,q+1}\,h^{\,q}.
\end{equation}
\end{rem}

\begin{proof}
Let $\tau_j=h\tilde{s}_j$ and $\tilde{s}_j\in [0,1]$. By \cref{lem:pullout}, we may drop outside conjugations, so it suffices to analyze
nested commutators of the occurrences $H(\tau_j)=V_{h\tilde s_j}$.
\begin{equation}
B_{h\tilde s}:=H(h\tilde s)=e^{\frac{i}{h}h\tilde s\,A}\,V\,e^{-\frac{i}{h}h\tilde s\,A}=V_{h\tilde s}.
\end{equation}
Unitary transformations preserve $L^2$ norms, so it suffices to analyze nested
commutators built only from the $B_{h\tilde s_j}$.

By the exact Egorov Theorem as in \cref{lem:egorov}, we have
\begin{equation}
V_A(h\tilde s)=e^{\frac{i}{h}h\tilde s A} \, V \, e^{-\frac{i}{h}h\tilde s A}
=\op_h\!\big(v(x-\tilde s p)\big)=\op_h(w_{\tilde s}),\quad
w_{\tilde s}(x,p)=v(x-p\tilde s)\in S(1).
\end{equation}
Because $(x,p)\mapsto(x-\tilde s p,p)$ preserves $S(1)$ for $0\le \tilde s\le1$, we have $w_{\tilde s}\in S(1)$ uniformly. This establishes the
case $p=1$ of \cref{eq:main-symbol}, which is the base case of $p=1$.

By \cref{lem:left-normed-cont}, any left-normed commutator with $n+1$ occurrences
(hence $n$ layers) satisfies
\begin{equation} \label{eq:comm_n+1}
\Comm_{n+1}(h\tilde s_1,\ldots,h\tilde s_{n+1})
= h^{\,n}\,\op_h(\tilde S_{n+1}),\qquad
\tilde S_{n+1}\in S(1)\ \ \text{in }(\tilde s_1,\ldots,\tilde s_{n})\in[0,1]^n.
\end{equation}
By \cref{lem:pullout}, we have
\begin{equation}
\operatorname{ad}_{H(t_{\tau_1})}\cdots\operatorname{ad}_{H(\tau_n)}\big(H(\tau_{n+1})\big)
= e^{\,iA\tau_{n+1}}\,\Comm_{n+1}(\tau_{1}-\tau_{n+1},\dots,\tau_{n}-\tau_{n+1})\,e^{-\,iA\tau_{n+1}},
\end{equation}
By \cref{eq:comm_n+1} together with \cref{lem:egorov}, we have 
\begin{equation}\label{eq:Hcommh} \operatorname{ad}_{H(t_{\tau_1})}\cdots\operatorname{ad}_{H(\tau_n)}\big(H(\tau_{n+1})\big)=h^{n} \op_h(g_{n+1}), 
\end{equation}
where $g_{n+1}\in S(1)$.

Next, write the given $C\!\big(H(\tau_1),\ldots,H(\tau_p)\big)\in \mathcal{C}_p$ as an iterated commutator
of $k\ge1$ blocks:
\begin{equation}
    C=\big[\;\cdots\big[\;
\underbrace{C_{n_1}}_{\text{$n_1$ occurrences}},
\underbrace{C_{n_2}}_{\text{$n_2$ occurrences}}\big]\ ,\ldots\ ,\,
\underbrace{C_{n_k}}_{\text{$n_k$ occurrences}}\big],
\end{equation}
where each $C_{n_j}$ is has $n_j$ occurrences of $H(\cdot)$ in the form of \cref{eq:Hcommh}, where the total number of occurrences of $H(t)$ satisfies
\begin{equation}
n_1+\cdots+n_k=p.
\end{equation}
so the total number of layers is $(n_1-1)+\cdots+(n_k-1)+(k-1)=p-1$.
Thus, by \cref{eq:Hcommh}, each $C_{n_j}$ contributes a factor of $h^{n_j}$. In addition, each outer commutator contributes one additional factor of $h$ and keeps the symbol in $S(1)$
by the Weyl commutator rule in \cref{lem:comm_rule}, and 
commuting the $k$ blocks produces $h^{k-1}$. In sum, we have 
\begin{equation}\label{eq:main-symbol}
C = h^{\,(n_1-1)+\cdots+(n_k-1)}\,h^{\,k-1}\,\op_h\!\big(\tilde{g}\big)
= h^{\,p-1}\,\op_h\!\big(\tilde{g}\big),
\end{equation}
for some symbol $\tilde{g}\in S(1)$ whose $S(1)$ seminorms are uniformly bounded in the
times $\tilde s_1,\ldots,\tilde s_p\in [0,1]$.
Finally, by the Calder\'on-Vaillancourt theorem (\cref{thm:CV}), we have the desired result.

\end{proof}

\begin{rem}[On attempting a Taylor expansion]
Although it seems that we can use the Taylor expansion of $H(t)$ at $t_j$,
\begin{equation}
H(t)=H(t_j)+H'(t_j)(t-t_j)+\cdots,
\end{equation}
and it seems this would also gain some order in $h$, there are two issues. Crucially, the leading first nonzero term is not of order $h^p$. For example, take three occurrences (a two layer commutator). With $t_\alpha\in  (t_j, t_{j+1})$ where $\alpha=1,2,3$.
Expanding around $t_j$ gives
\begin{equation}
H(t_\alpha)=H(t_j)+(t_\alpha-t_j)H'(t_j)+\mathcal{O}(h^2),
\qquad t_\alpha-t_j=\mathcal{O}(h).
\end{equation}
Setting $H:=H(t_j)$ and $H':=H'(t_j)$, one finds
\begin{equation}
\begin{aligned}
[H(t_1),[H(t_2),H(t_3)]]
&=[H+(t_1-t_j)H',\,(t_3-t_2)[H,H']]+\mathcal{O}(h^2) \\
&=(t_3-t_2)\Big([H,[H,H']]+(t_1-t_j)[H',[H,H']]\Big)+\mathcal{O}(h^2), \\
\end{aligned}
\end{equation}
which implies
\begin{equation}
[H(t_1),[H(t_2),H(t_3)]]=\mathcal{O}(h).
\end{equation}
Thus the Taylor expansion yields only a single factor of $h$. In particular, for two layers one does not reach the $h^2$ scaling required by \cref{thm:main-L2}. In addition, this approach will introduce polynomial dependence on $H'$ (and higher derivatives), so the resulting bounds are no longer uniform in the norm of the Hamiltonian's time derivative,
unlike the clean estimates provided by \cref{thm:main-L2}. Note that an important feature of the quantum Magnus algorithms is to have logarithmic dependence on the Hamiltonian's time derivative, instead of a polynomial dependence. Therefore, although Taylor expansion suggests a partial gain in order, it does not yield superconvergence.
\end{rem}

\begin{thm}[Global $2p$-order for the $p$-th order Magnus algorithm in the interaction picture]\label{thm:global-2p}
Let $A=-\tfrac12\Delta$ and $B=V(x)$ with $V\in S(1)$. Let $H(t)$ be the interaction-picture Hamiltonian given by 
\begin{equation}
H(t)=e^{iAt}\,B\,e^{-iAt}.
\end{equation}
Let $U(T,0)$ denote the exact evolution operator, and let $U_p(T,0)$ denote the global $p$-th order Magnus approximation with step size $h=T/L$, defined as in \cref{eq:U} and \cref{eq:U_p}, respectively.
Then there exists a constant $C_{V,p}>0$ such that for all $0<h\le 1$, we have \cref{eq:def_alpha_comm_q}
\begin{equation}\label{eq:global-2p-bound}
\bigl\|U(T,0)-U_p(T,0)\bigr\|_{L^2(\RR^d)\to L^2(\RR)} \le\ C_{V,p}\,T\,h^{2p}.
\end{equation}
The constant $C_{V,p}$ depends only on $p$, the dimension $d$, and finitely many $L^\infty$ bounds of spatial derivatives of $V$ (i.e.\ on $S(1)$ seminorms of $V$).
\end{thm}

\begin{proof}
By the Magnus local truncation error bound in \cref{thm:magnus_trunction_finite_comm}, it is sufficient to consider $\alpha_{\mathrm{comm},q}$. In particular, we have
\begin{equation} \label{eq:gradeq}
\alpha_{\mathrm{comm},q} = \sup_{\tau_1, \ldots, \tau_q \in [t_j,t_{j+1}]} \max_{C \in \mathcal{C}_q} \big\| C\big(H(\tau_1), \ldots, H(\tau_q)\big) \big\| \leq  C_{V,q}\, h^{\,q-1},
 \end{equation}
 by \cref{thm:main-L2}.
Substituting \cref{eq:gradeq} into \cref{eq: LTE bound} yields
\begin{equation}
\norm{U(t_{j+1},t_j) - U_p(t_{j+1},t_j)}
\ \le\ C\,\sum_{q=p+1}^{p^2+2p} C_{V,q}\,h^{\,q-1}\, h^{\,q}
\ =\ C\,\sum_{q=p+1}^{p^2+2p} C_{V,q}\, h^{\,2q-1}.
\end{equation}
Because $h\in(0,1]$ and the index set is finite, there exists a constant $\tilde{C}_{V,p}>0$ such that
\begin{equation}
\sum_{q=p+1}^{p^2+2p} C_{V,q}\, h^{\,2q-1}  \leq \sum_{q=p}^{p^2+2p-1} C_{V,q+1} \,h^{\,2q+1}\ \le\ \tilde{C}_{V,p}\,h^{\,2p+1},
\end{equation}
which estimates the local truncation error.
Absorbing $C$ into $\tilde{C}_{V,p}$ and since $U$ and $U_p$ are both unitary, summing the local errors over $L$ time steps then yields the global error bound in \cref{eq:global-2p-bound}, as desired. 
\end{proof}

\section{Numerical Results}
In this section, we present numerical evidence supporting the theoretical $h$-scaling bounds for time-labeled nested commutators in the interaction picture. We consider
\begin{equation} \label{eq:setup_ABV}
A=-\tfrac12\Delta,\qquad B=V(x),\qquad V(x)=\cos x,
\end{equation}
on $[-\pi,\pi)$ with periodic boundary conditions in one spatial dimension, and
\begin{equation}
H(t)=e^{iAt} B e^{-iAt}.
\end{equation}
Theorem~\ref{thm:main-L2} predicts that any grade-$p$ commutator of $\{H(t)\}$, which has $p-1$ layers, obeys
\begin{equation}
\|C(H(\tau_1),\dots,H(\tau_p))\|_{L^2\to L^2}\le C_{V,p} h^{p-1}
\end{equation}
uniformly in the time labels.

The Laplacian is discretized by the second-order central finite difference scheme on a uniform grid with
\begin{equation}
N\in\{256,512,1024,2048\}.
\end{equation}
The potential is applied diagonally, and all operator sizes are reported in the matrix spectral norm.
We sample the semiclassical parameter (equivalently the time step size) as
\begin{equation}
h\in\{1,2^{-1},2^{-2},2^{-3},2^{-4},2^{-5}\}\quad\text{(that is, }h=2^{-k},\ k=0,\dots,5\text{)}.
\end{equation}
For each $h$ we use the time-label set 
\begin{equation}
\{h,\ h/2,\ h/4,\ h/8,\ h/16,\ h/32,\ h/64 \}\subset[0,h].
\end{equation}
For the three-layer case, for every triple $(\tau,s,\sigma)$ drawn from this set, we evaluate the commutator $[H(\tau),[H(s),H(\sigma)]]$, and for each choice of $(h,N)$, we record the largest spectral norm observed across all triples. For the four-layer case, the same procedure is applied to $[H(\tau),[H(s),[H(\sigma),H(\rho)]]]$ over all quadruples $(\tau,s,\sigma,\rho)$ from the same set. Thus, for each $(h,N)$, the plotted value is the maximum spectral norm taken over all corresponding time labels, ensuring a consistent and reproducible procedure across different parameters.

\begin{figure}[!ht]
    \centering
    \begin{subfigure}[t]{0.47\linewidth}
        \centering
        \includegraphics[width=\linewidth]{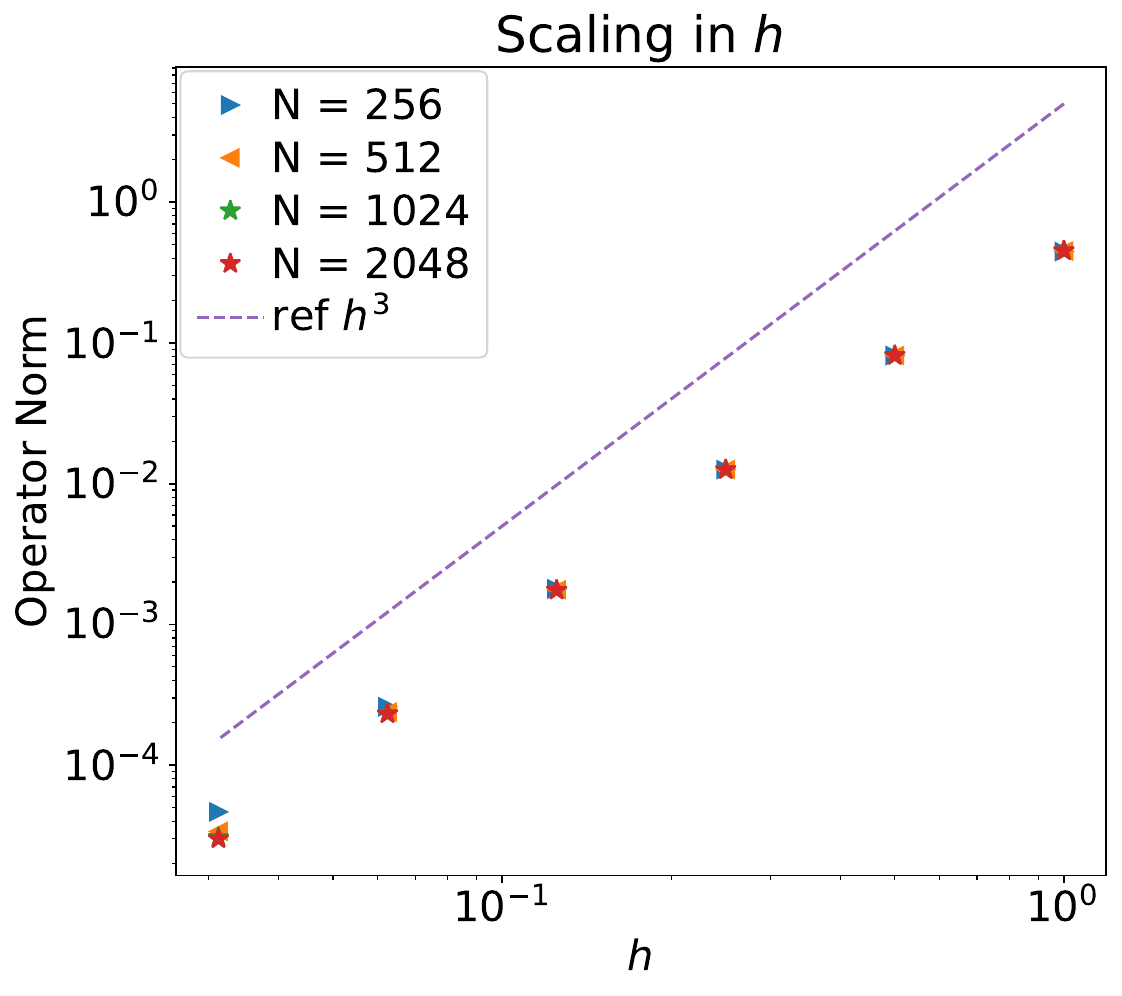}
        \caption{Scaling in $h$ for the three-layer commutator}
        \label{fig:comm3}
    \end{subfigure}
    \hfill
    \begin{subfigure}[t]{0.47\linewidth}
        \centering
        \includegraphics[width=\linewidth]{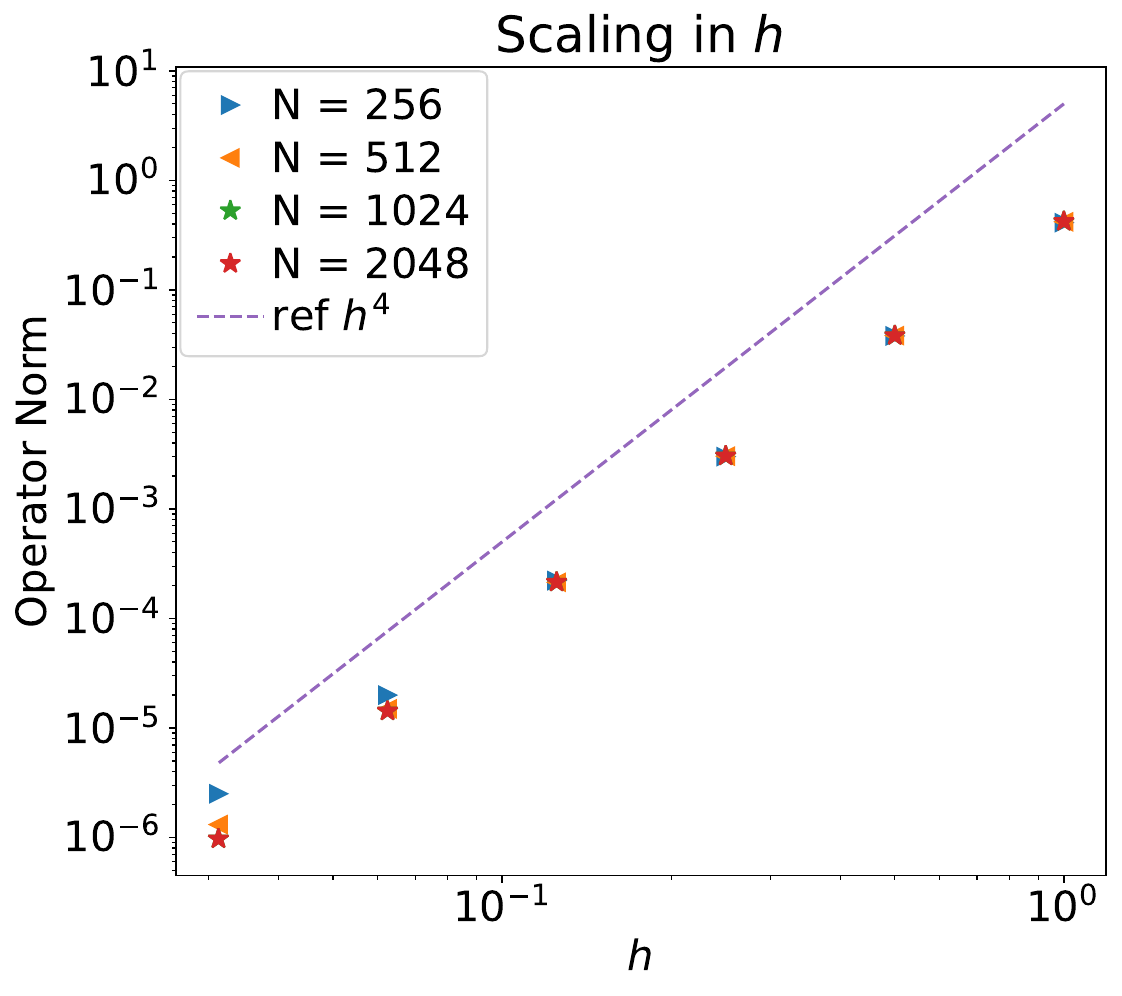}
        \caption{Scaling in $h$ for the four-layer commutator}

        \label{fig:comm4}
    \end{subfigure}
    \caption{Log-log plots of the largest sampled spectral norms of nested commutators versus $h\in\{1,2^{-1},2^{-2},2^{-3},2^{-4},2^{-5}\}$. Left: three layers, where for each $h$ and $N\in\{256,512,1024,2048\}$ we evaluate the maximum of the spectral norm of $[H(\tau),[H(s),H(\sigma)]]$ over $(\tau,s,\sigma) \in[0,h]$. Right: four layers, defined analogously for $[H(\tau),[H(s),[H(\sigma),H(\rho)]]]$ over $(\tau,s,\sigma,\rho)$. Here $H(t)=e^{iAt} B e^{-iAt}$ with setup in \cref{eq:setup_ABV}. Axes are labeled $h$ (horizontal, log scale) and operator norm (vertical, log scale). The dashed lines of slopes 3 and 4 mark the expected $O(h^3)$ and $O(h^4)$ behavior. The near-coincident curves across $N$ indicate uniformity in $N$.}
    \label{fig:comm_scaling}
\end{figure}

In \cref{fig:comm_scaling}, the horizontal axis is $h$ (log scale) and the vertical axis is operator norm (log scale). Each colored curve corresponds to one grid size $N\in\{256,512,1024,2048\}$. A dashed reference line with slope 3 is shown in the three-layer subplot and a dashed reference line with slope 4 in the four-layer subplot to indicate the expected $\mathcal{O}(h^3)$ and $\mathcal{O}(h^4)$ rates from \cref{thm:main-L2} with $p=4$ and $p=5$. In both subplots, the curves for different $N$ nearly coincide at each $h$, demonstrating that the measured $h$-rates are uniform with respect to the spatial resolution on the tested grids. 

Having established the $h$-scaling of time-labeled nested commutators, we next demonstrate the behavior of the Magnus algorithm itself on the same model to connect the bounds with a complete Magnus approximation. We keep the spatial discretization and norms as above and fix the grid at $N=128.$
In the following, we write the time step as $\Delta t$ and report a single-step error (local truncation error) over 
$[t_j,t_j+\Delta t]$.

With $A=-\tfrac12\Delta$ and $B=V(x)$, the exact one-step propagator we compare against is the interaction picture unitary given by \cref{eq:U_IP}. Note that directly simulating this on a quantum device using any time-independent Hamiltonian simulation algorithm has a cost scale at least linearly in the operator norm of $A$ (equivalently $\mathrm{poly}(N)$), while our new quantum algorithm is $\mathrm{polylog}(N)$.
For the Magnus approximations, we use the same truncated per-step unitaries as in Algorithm~1,
\begin{equation}
U_p(\Delta t)=\exp\!\big(\Omega_{(p)}(\Delta t)\big),\qquad
\Omega_{(p)}(\Delta t)=\sum_{n=1}^{p}\Omega_n(\Delta t),
\end{equation}
and in the experiments below we instantiate $p=1,2$, i.e.
\begin{equation}
U_1(\Delta t)=\exp\!\big(\Omega_1(\Delta t)\big),\qquad
U_2(\Delta t)=\exp\!\big(\Omega_1(\Delta t)+\Omega_2(\Delta t)\big).
\end{equation}
The reported error is the spectral norm
\begin{equation}
\bigl\|U_p(\Delta t)-U_{\mathrm{exact}}(\Delta t)\bigr\|_2,\qquad p\in\{1,2\}.
\end{equation}
Time integrals in $\Omega_n$ are evaluated by Gauss-Legendre (GL) quadrature on $[0,\Delta t]$ in the eigenbasis of $A$. We denote by $M$ the number of GL nodes per time integral and use
\begin{equation}
M_1=512\ \text{for }\Omega_1,\qquad
M_2=256\ \text{(per layer) for }\Omega_2,
\end{equation}
with lower-triangular ordering to enforce $t_2\le t_1$. Note that in the quantum Magnus algorithm, the computational cost depends logarithmically on the number of quadrature points
$M$. We vary the time step over
\begin{equation}
\Delta t\in\{0.8,\,0.4,\,0.2,\,0.1\}.
\end{equation}

\begin{figure}[!h]
    \centering
    \begin{subfigure}[t]{0.49\linewidth}
        \centering
        \includegraphics[width=\linewidth]{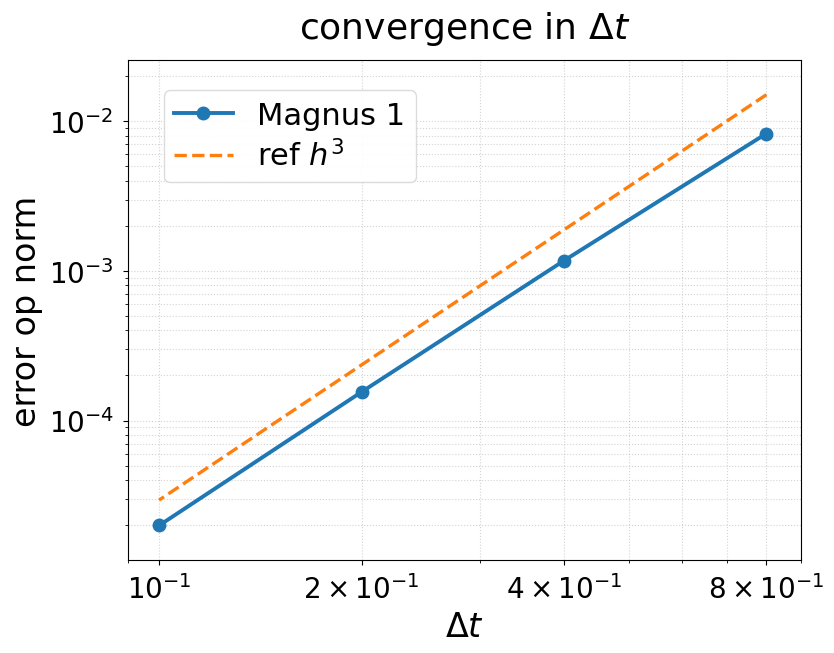}
        \caption{Convergence in $\Delta t$ for $U_1=\exp(\Omega_1)$}
        \label{fig:Mag1}
    \end{subfigure}
    \hfill
    \begin{subfigure}[t]{0.49\linewidth}
        \centering
        \includegraphics[width=\linewidth]{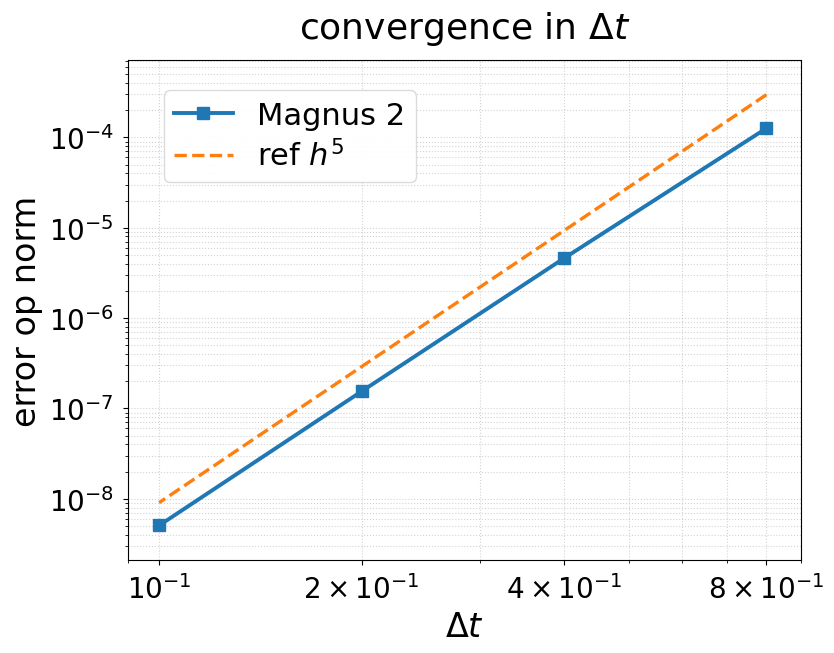}
        \caption{Convergence in $\Delta t$ for $U_2=\exp(\Omega_1+\Omega_2)$}
        \label{fig:Mag2}
    \end{subfigure}
    \caption{Log-log plots of the one-step unitary error
    $\|U_p(\Delta t)-U_{\mathrm{exact}}(\Delta t)\|_2$ versus $\Delta t\in\{0.8,0.4,0.2,0.1\}$ for $p=1$ (left) and $p=2$ (right).
    Spatial grid: $N=128$. Gauss-Legendre quadrature for Magnus time integrals: $M_1=512$ for $\Omega_1$, $M_2=256$ per layer for $\Omega_2$. Potential: $V(x)=\tfrac12\cos x$.
    The dashed reference lines of slopes $3$ (left) and $5$ (right) indicate the expected $\Or(\Delta t^3)$ and $\Or(\Delta t^5)$ local truncation behavior, respectively. In both subplots the measured slopes align with the references. Axes are $\Delta t$ (horizontal, log scale) and the operator norm of error (vertical, log scale).}
    \label{fig:Magnus}
\end{figure}

According to \cref{thm:global-2p}, $U_1$ yields a global error of order $\Or(\Delta t^2)$ and a local truncation error of order $\Or(\Delta t^3)$, while $U_2$ yields a global error of order $\Or(\Delta t^4)$ and a local truncation error of order $\Or(\Delta t^5)$. Since our plots measure the single-step error in $[0,\Delta t]$, the expected slopes in the log-log scaling are $3$ for $U_1$ as in the left subplot of \cref{fig:Magnus} and $5$ for $U_2$ as in the right subplot of \cref{fig:Magnus}. The numerical curves in \cref{fig:Magnus} align precisely with these reference lines, confirming the predicted one-step convergence rates and providing an algorithmic complement to the commutator scaling established earlier.

\section{Conclusion and Remarks}
In this work, we show that high-order quantum Magnus schemes for interaction picture Schr\"odinger dynamics exhibit a uniform order-doubling phenomenon. For $A=-\tfrac12\Delta$ and $B=V(x)$ with $V\in S(1)$, the global $p$-th order method achieves accuracy $\Or(h^{2p})$ in operator norm with constants depending only on finitely many derivatives of $V$ and independent of the kinetic norm and of the spatial resolution. The analysis reorganizes time labels so that all nested commutators of $H_I(t)=e^{iAt}Be^{-iAt}$ reduce to commutators among shifted observables $B_u=e^{iAu}Be^{-iAu}$ with a fixed innermost entry. The exact Egorov Theorem for the quadratic kinetic energy transports symbols as $B_u=\op_h\!\big(v(x-u p)\big)$, and an iterative Poisson bracket expansion yields one power of $h$ per commutator layer. The Calder\'on-Vaillancourt theorem then converts uniform symbol bounds into uniform operator norm bounds. Together with the circuit constructions of~\cite{FangLiuZhu2025}, whose per-step gate complexity scales $\mathrm{polylog}(N)$, the $2p$ order estimate leads to a total gate complexity that remains $\mathrm{polylog}(N)$ at fixed target accuracy.

A natural next step is to transfer the continuum analysis to fully discrete evolutions. As shown in earlier discrete work~\cite{Borns-WeilFangZhang2025,BornsWeilFang2022}, one has the discrete to continuous operator comparison
\begin{equation}\label{eq:disc-cont}
\big\|\op_{N}(a)\big\|_{H_N\to H_N}\;\le\; \big\|\op_h(\tilde a)\big\|_{L^2(\RR^d)\to L^2(\RR^d)} ,
\end{equation}
where $\op_{N}(a)$ denotes the discrete operator obtained from discretization on $N$ points and $\op_h(\tilde a)$ is the corresponding semiclassical pseudodifferential operator acting on $L^2(\RR^d)$.
This inequality provides a clear route: first transfer the uniform commutator estimates derived here to the discrete setting via \cref{eq:disc-cont}, and then propagate them through the discrete Magnus error representation. We anticipate technical challenges that are specific to the discrete version, such as developing a two-parameter symbol framework that simultaneously tracks spatial resolution and time step. A further challenge is to handle periodic extensions and aliasing effects in formulations based on torus lifts and related discretizations. Addressing these issues fits naturally with the methods introduced here and in~\cite{BornsWeilFang2022,Borns-WeilFangZhang2025}, and we will pursue them in future work.

There are further directions that complement the discrete extension. One is to relax smoothness by replacing $V\in S(1)$ with other regularity classes or piecewise smooth potentials and to quantify how the number of required derivatives depends on the Magnus order. Another is to treat kinetic energies beyond the quadratic case, including magnetic flows, where the Egorov Theorem is only approximate. In this setting, one possible approach is to track the growth of symbols under repeated conjugations and control the resulting losses up to Ehrenfest time intervals. Combined with the commutator hierarchy, this may still yield effective order doubling on time scales of physical relevance. 
A third is to study robustness under a smaller number of quadrature points, including sparse or randomized time sampling, and its interaction with the cancellation structure in the nested commutators for other interesting applications. Each of these topics presents concrete mathematical and algorithmic questions that warrant separate investigation, which we will leave for future work.

In summary, we provide a uniform, symbol calculus explanation for superconvergence of high-order quantum Magnus algorithms in the interaction picture and a concrete bridge to discrete simulations through \eqref{eq:disc-cont} and the framework of~\cite{Borns-WeilFangZhang2025,BornsWeilFang2022}. We expect these ideas to inform rigorous error analysis and practical algorithm design for a broad range of time-dependent quantum simulations.

\section*{Acknowledgments}
This work is supported by the National Science Foundation CAREER Award DMS-2438074 (D.F. and J.Z.), and the U.S. Department of Energy, Office of Science, Accelerated Research in Quantum Computing Centers, Quantum Utility through Advanced Computational Quantum Algorithms, grant no. DE-SC0025572.

\bibliographystyle{unsrt}
\bibliography{magnus-p}

\begin{thebibliography}{10}

\bibitem{ChildsSuTranEtAl2020}
Andrew~M. Childs, Yuan Su, Minh~C. Tran, Nathan Wiebe, and Shuchen Zhu.
\newblock Theory of trotter error with commutator scaling.
\newblock {\em Phys. Rev. X}, 11:011020, 2021.

\bibitem{ChildsSu2019}
Andrew~M Childs and Yuan Su.
\newblock Nearly optimal lattice simulation by product formulas.
\newblock {\em Phys. Rev. Lett.}, 123(5):050503, 2019.

\bibitem{WatkinsWiebeRoggeroLee2022}
Jacob Watkins, Nathan Wiebe, Alessandro Roggero, and Dean Lee.
\newblock Time-dependent hamiltonian simulation using discrete-clock constructions.
\newblock {\em PRX Quantum}, 5(4):040316, 2024.

\bibitem{ZhukRobertsonBravyi2023}
Sergiy Zhuk, Niall Robertson, and Sergey Bravyi.
\newblock Trotter error bounds and dynamic multi-product formulas for hamiltonian simulation, 2023.

\bibitem{AftabAnTrivisa2024}
Junaid Aftab, Dong An, and Konstantina Trivisa.
\newblock Multi-product hamiltonian simulation with explicit commutator scaling, 2024.
\newblock arXiv:2403.08922.

\bibitem{Watson2025}
James~D. Watson.
\newblock Randomly compiled quantum simulation with exponentially reduced circuit depths, 2025.

\bibitem{WatsonWatkins2024}
James~D. Watson and Jacob Watkins.
\newblock Exponentially reduced circuit depths using trotter error mitigation.
\newblock {\em PRX Quantum}, 6:030325, Aug 2025.

\bibitem{Somma2015}
Rolando~D. Somma.
\newblock Quantum simulations of one dimensional quantum systems, 2015.
\newblock arXiv:1503.06319.

\bibitem{SahinogluSomma2020}
Burak Şahinoğlu and Rolando~D. Somma.
\newblock Hamiltonian simulation in the low-energy subspace.
\newblock {\em npj Quantum Information}, 7(1), July 2021.

\bibitem{AnFangLin2021}
Dong An, Di~Fang, and Lin Lin.
\newblock Time-dependent unbounded {H}amiltonian simulation with vector norm scaling.
\newblock {\em {Quantum}}, 5:459, may 2021.

\bibitem{SuHuangCampbell2021}
Yuan Su, Hsin-Yuan Huang, and Earl~T. Campbell.
\newblock Nearly tight {T}rotterization of interacting electrons.
\newblock {\em {Quantum}}, 5:495, July 2021.

\bibitem{ZhaoZhouShawEtAk2021}
Qi~Zhao, You Zhou, Alexander~F. Shaw, Tongyang Li, and Andrew~M. Childs.
\newblock Hamiltonian simulation with random inputs.
\newblock {\em Physical Review Letters}, 129(27), December 2022.

\bibitem{ChildsLengEtAl2022}
Andrew~M. Childs, Jiaqi Leng, Tongyang Li, Jin-Peng Liu, and Chenyi Zhang.
\newblock Quantum simulation of real-space dynamics.
\newblock {\em Quantum}, 6:860, November 2022.

\bibitem{FangTres2023}
Di~Fang and Albert Tres~Vilanova.
\newblock Observable error bounds of the time-splitting scheme for quantum-classical molecular dynamics.
\newblock {\em SIAM J. Numer. Anal.}, 61(1):26--44, 2023.

\bibitem{BornsWeilFang2022}
Yonah Borns-Weil and Di~Fang.
\newblock Uniform observable error bounds of trotter formulae for the semiclassical schrödinger equation.
\newblock {\em Multiscale Modeling \& Simulation}, 23(1):255--277, 2025.

\bibitem{HuangTongFangSu2023}
Hsin-Yuan Huang, Yu~Tong, Di~Fang, and Yuan Su.
\newblock Learning many-body hamiltonians with heisenberg-limited scaling.
\newblock {\em Phys. Rev. Lett.}, 130:200403, May 2023.

\bibitem{ZengSunJiangZhao2022}
Pei Zeng, Jinzhao Sun, Liang Jiang, and Qi~Zhao.
\newblock Simple and high-precision hamiltonian simulation by compensating trotter error with linear combination of unitary operations, 2022.
\newblock arXiv:22212.04566.

\bibitem{GongZhouLi2023}
Weiyuan Gong, Shuo Zhou, and Tongyang Li.
\newblock Complexity of {D}igital {Q}uantum {S}imulation in the {L}ow-{E}nergy {S}ubspace: {A}pplications and a {L}ower {B}ound.
\newblock {\em {Quantum}}, 8:1409, July 2024.

\bibitem{LowSuTongTran2023}
Guang~Hao Low, Yuan Su, Yu~Tong, and Minh~C. Tran.
\newblock Complexity of implementing trotter steps.
\newblock {\em PRX Quantum}, 4:020323, May 2023.

\bibitem{ZhaoZhouChilds2024}
Qi~Zhao, You Zhou, and Andrew~M Childs.
\newblock Entanglement accelerates quantum simulation.
\newblock {\em Nature Physics}, pages 1--8, 2025.

\bibitem{YuXuZhao2024}
Wenjun Yu, Jue Xu, and Qi~Zhao.
\newblock Observable-driven speed-ups in quantum simulations, 2024.

\bibitem{ChenXuZhaoYuan2024}
Boyang Chen, Jue Xu, Qi~Zhao, and Xiao Yuan.
\newblock Error interference in quantum simulation, 2024.
\newblock arXiv:2411.03255.

\bibitem{BerryChildsCleveEtAl2015}
D.~W. Berry, A.~M. Childs, R.~Cleve, R.~Kothari, and R.~D. Somma.
\newblock Simulating {Hamiltonian} dynamics with a truncated {Taylor} series.
\newblock {\em Phys. Rev. Lett.}, 114:090502, 2015.

\bibitem{KieferovaSchererBerry2019}
Mária Kieferová, Artur Scherer, and Dominic~W. Berry.
\newblock Simulating the dynamics of time-dependent hamiltonians with a truncated dyson series.
\newblock {\em Physical Review A}, 99(4), Apr 2019.

\bibitem{LowWiebe2019}
G.~H. Low and N.~Wiebe.
\newblock Hamiltonian simulation in the interaction picture.
\newblock 2019.
\newblock arXiv:1805.00675.

\bibitem{BerryChildsSuEtAl2020}
D.~W. Berry, A.~M. Childs, Y.~Su, X.~Wang, and N.~Wiebe.
\newblock Time-dependent {H}amiltonian simulation with $l^{1}$-norm scaling.
\newblock {\em Quantum}, 4:254, 2020.

\bibitem{LowChuang2017}
Guang~Hao Low and Isaac~L. Chuang.
\newblock Optimal {H}amiltonian simulation by quantum signal processing.
\newblock {\em Phys. Rev. Lett.}, 118:010501, 2017.

\bibitem{LowChuang2019}
Guang~Hao Low and Isaac~L. Chuang.
\newblock Hamiltonian simulation by qubitization.
\newblock {\em Quantum}, 3:163, Jul 2019.

\bibitem{GilyenSuLowEtAl2019}
Andr{\'a}s Gily{\'e}n, Yuan Su, Guang~Hao Low, and Nathan Wiebe.
\newblock Quantum singular value transformation and beyond: exponential improvements for quantum matrix arithmetics.
\newblock In {\em Proceedings of the 51st Annual ACM SIGACT Symposium on Theory of Computing}, pages 193--204, 2019.

\bibitem{ZlokapaSomma2024}
Alexander Zlokapa and Rolando~D. Somma.
\newblock Hamiltonian simulation for low-energy states with optimal time dependence.
\newblock {\em Quantum}, 8:1449, August 2024.

\bibitem{ZhengLengLiuWu2024}
Yufan Zheng, Jiaqi Leng, Yizhou Liu, and Xiaodi Wu.
\newblock On the computational complexity of schr\"odinger operators, 2024.

\bibitem{AnFangLin2022}
Dong An, Di~Fang, and Lin Lin.
\newblock Time-dependent hamiltonian simulation of highly oscillatory dynamics and superconvergence for schr{\"o}dinger equation.
\newblock {\em Quantum}, 6:690, 2022.

\bibitem{FangLiuSarkar2025}
Di~Fang, Diyi Liu, and Rahul Sarkar.
\newblock Time-dependent hamiltonian simulation via magnus expansion: Algorithm and superconvergence.
\newblock {\em Communications in Mathematical Physics}, 406(6):128, 2025.

\bibitem{FangLiuZhu2025}
Di~Fang, Diyi Liu, and Shuchen Zhu.
\newblock High-order magnus expansion for hamiltonian simulation, 2025.
\newblock arXiv:2509.06054.

\bibitem{Borns-WeilFangZhang2025}
Yonah Borns-Weil, Di~Fang, and Jiaqi Zhang.
\newblock {Discrete Superconvergence Analysis for Quantum Magnus Algorithms of Unbounded Hamiltonian Simulation}, 2025.
\newblock arXiv:2502.20255.

\bibitem{HochbruckLubich2003}
Marlis Hochbruck and Christian Lubich.
\newblock On {M}agnus integrators for time-dependent {S}chr\"{o}dinger equations.
\newblock {\em SIAM J. Numer. Anal.}, 41(3):945--963, 2003.

\bibitem{Thalhammer2006}
Mechthild Thalhammer.
\newblock A fourth-order commutator-free exponential integrator for nonautonomous differential equations.
\newblock {\em SIAM Journal on Numerical Analysis}, 44(2):851--864, 2006.

\bibitem{BlanesCasasThalhammer2017}
Sergio Blanes, Fernando Casas, and Mechthild Thalhammer.
\newblock High-order commutator-free quasi-magnus exponential integrators for non-autonomous linear evolution equations.
\newblock {\em Computer Physics Communications}, 220:243--262, 2017.

\bibitem{BlanesCasasOteoRos2009}
S.~Blanes, F.~Casas, J.~A. Oteo, and J.~Ros.
\newblock The {M}agnus expansion and some of its applications.
\newblock {\em Phys. Rep.}, 470(5-6):151--238, 2009.

\bibitem{SharmaTran2024}
Kunal Sharma and Minh~C. Tran.
\newblock Hamiltonian simulation in the interaction picture using the magnus expansion, 2024.
\newblock arXiv:2404.02966.

\bibitem{BosseChildsEtAl2024}
Jan~Lukas Bosse, Andrew~M Childs, Charles Derby, Filippo~Maria Gambetta, Ashley Montanaro, and Raul~A Santos.
\newblock Efficient and practical hamiltonian simulation from time-dependent product formulas.
\newblock {\em Nature Communications}, 16(1):2673, 2025.

\bibitem{CasaresZiniArrazola2024}
Pablo Antonio~Moreno Casares, Modjtaba~Shokrian Zini, and Juan~Miguel Arrazola.
\newblock Quantum simulation of time-dependent hamiltonians via commutator-free quasi-magnus operators, 2024.
\newblock arXiv:2403.13889.

\bibitem{ApelCubittOnorati2025}
Harriet Apel, Toby Cubitt, and Emilio Onorati.
\newblock A sharper magnus expansion bound woven in binary branches, 2025.
\newblock arXiv:2509.18312.

\bibitem{zworski2022semiclassical}
Maciej Zworski.
\newblock {\em Semiclassical Analysis}.
\newblock American Mathematical Society, 2012.

\bibitem{Martinez2002}
Andr{\'e} Martinez.
\newblock {\em {An introduction to semiclassical and microlocal analysis}}, volume 994.
\newblock Springer, 2002.

\bibitem{BlanesCasas2017book}
Sergio Blanes and Fernando Casas.
\newblock {\em A concise introduction to geometric numerical integration}.
\newblock CRC press, 2017.

\bibitem{HairerHochbruckIserlesLubich2006}
Ernst Hairer, Marlis Hochbruck, Arieh Iserles, and Christian Lubich.
\newblock Geometric numerical integration.
\newblock {\em Oberwolfach Reports}, 3(1):805--882, 2006.

\end{thebibliography}
\end{document}